\documentclass[12pt]{amsart}

\usepackage{fullpage}
\usepackage[british]{babel}  
\usepackage{geometry}    
\usepackage{mathrsfs}
\usepackage{bbm}					
\usepackage[parfill]{parskip} 
 \usepackage{graphicx}				
\usepackage{breqn}
\usepackage{bm}	
\usepackage{amsmath,amssymb,amsthm,eucal,verbatim}
\usepackage{hyperref,amsbsy}
\usepackage{graphicx}
\usepackage{epstopdf}

\usepackage{enumitem}

\usepackage{mathtools}

\advance \topmargin by-\headheight
\advance \topmargin by-\headsep
\evensidemargin \oddsidemargin
\newtheorem{theorem}{Theorem}[section]
\newtheorem{lem}{Lemma}[section]
\newtheorem*{thm*}{Theorem}
\newtheorem*{lem*}{Lemma}
\newtheorem*{cor*}{Corollary}

\newtheorem*{rem*}{Remark}

\theoremstyle{definition}

\newtheorem*{thm A}{Theorem A}

\newtheorem*{hyp}{Hypothesis $\mathscr{H}_{\delta}$}

\newtheorem{conj}{Conjecture}

\newcommand{\be}{\begin{equation}}
\newcommand{\ee}{\end{equation}}
\renewcommand{\a}{\alpha}

\renewcommand{\b}{\beta}
\renewcommand{\r}{\rho }
\newcommand{\g}{\gamma }
\newcommand{\G}{\Gamma }

\newcommand{\e}{\epsilon}
\renewcommand{\l}{\lambda}

\renewcommand{\d}{\delta}   
\newcommand{\sgn}{\operatorname{sgn}}

\makeatletter
\newcommand*{\bigs}[1]{{\hbox{$\left#1\vbox to10\p@{}\right.\n@space$}}}
\makeatother

\makeatletter
\newcommand*{\bigss}[1]{{\hbox{$\left#1\vbox to11\p@{}\right.\n@space$}}}
\makeatother

\newcommand{\sdfrac}[2]{\mbox{\small$\displaystyle\frac{#1}{#2}$}}

\title[The Uniform Distribution Modulo One of Certain Subsequences of Ordinates]
{The Uniform Distribution Modulo One of  Certain Subsequences of Ordinates of Zeros of the Zeta Function}

\author{Fatma \c{C}\.{i}\c{c}ek}
\address{Department of Mathematics and Statistics\\ University of Northern British Columbia\\ Prince George, BC, Canada}
\email{cicek@unbc.ca}

\author{Steven M. Gonek}
\address{Department of Mathematics\\ University of Rochester\\ Rochester, NY, USA}
\email{gonek@math.rochester.edu}

\calclayout

\begin{document}

\maketitle

\begin{abstract}

On the assumption of the Riemann hypothesis and a spacing hypothesis for the nontrivial zeros $\frac12+i\gamma$ of the Riemann zeta function, we show that the sequence 
\[
\Gamma_{[a, b]} =\Bigg\{ \gamma :  \gamma>0 \quad \mbox{and} \quad  
\frac{ \log\big(| \zeta^{(m_\g)} (\frac12+ i\g) | / (\log{\g} )^{m_\g}\big)}{\sqrt{\frac12\log\log \g}}  \in [a, b] \Bigg\},
\]
where the $\g$ are arranged in increasing order,  is uniformly distributed modulo one. Here $a$ and $b$ are real numbers with $a<b$, and $m_\g$ denotes the multiplicity of the zero $\frac12+i\g$. The same result holds when the $\g$'s are restricted to be the ordinates of  simple zeros. With   an extra hypothesis, we are also able to show an equidistribution result for   the scaled numbers $\gamma (\log T)/2\pi$ with $\g\in \Gamma_{[a, b]}$
and $0<\g\leq T$.
\end{abstract}


\section{Introduction}

It is well known that the  positive ordinates $\g$ of the  nontrivial zeros $\r=\b+i\g$ of the Riemann zeta function, when arranged in increasing order,  are uniformly distributed modulo one. 
This was proved by Rademacher~\cite{Rademacher} in the 1950s under the assumption of the Riemann hypothesis. Elliott~\cite{Elliott} later pointed out that this could be shown unconditionally. 
Our goal in this paper is to prove that if the Riemann hypothesis holds and a plausible hypothesis about the spacing of the $\g$'s is true, then the $\g$ are also uniformly distributed modulo one when we restrict to certain subsequences.

Throughout we assume the Riemann hypothesis so that every nontrivial zero of the zeta function has the form $\r=\frac12+i\g$. 
Then $N(T)$,  the number of ordinates $\g$ in the interval $(0, T]$ is given by  
\be\label{zero count}
N(T) =\frac{T}{2\pi} \log \frac{T}{2\pi}-\frac{T}{2\pi}+O\Big( \frac{\log T}{\log\log T}\Big)
\ee
(see Titchmarsh~\cite{T}, Ch. 14). Note that unconditionally, the error term is $O(\log T)$.  
We also assume that  for some $0<\d \leq 1$ the following spacing hypothesis holds for the zeros.
\begin{hyp} \label{hypothesis}
\emph{Let  $\g^+$ be the next larger ordinate of  a zero of the zeta function after the ordinate $\g$ with the understanding that 
     $\g^+=\g$ if and only if  $\frac12+i\g$ is a multiple zero.   Then there exists a positive constant $M$ such that,  uniformly for  $0<\l<1$, we have}
$$
 \limsup_{T\to\infty}\frac{1}{N(T)}\#\Big\{0 < \g \leq T:0\leq \g^+-\g \leq \frac{\l}{\log T}\Big\} \leq M \l^\d  .
$$
\end{hyp}

Hypothesis $\mathscr{H}_{\delta}$ is  credible because
Hypothesis $\mathscr{H}_{1}$ follows from Montgomery's pair correlation conjecture  
which, in turn, implies Hypothesis $\mathscr{H}_{\delta}$ for every $\delta \in (0, 1]$. Notice also 
that Hypothesis $\mathscr{H}_{\delta}$ implies that all but $o(N(T))$ of the zeros are simple, a fact we shall use later.

 Let $m_\g$ denote the multiplicity of the zero $\r=\frac12+i\g$, and let $a<b$ be real numbers. The sequences we wish to consider are
 $$
 \G_{[a, b] }= \Bigg\{  \g>0 \, : \,  \frac{ \log(\big| \zeta^{(m_\g)} (\frac12+ i\g) \big| / (\log{\g} )^{m_\g} )}{\sqrt{\frac12\log\log \g}}  \in [a, b] \Bigg\}
$$
and 
 $$
{ \G}^{*}_{[a, b] }= \Bigg\{  \g>0 \, : \, m_\g=1 \; \hbox{and}\; \; \frac{ \log(\big| \zeta^{'} (\frac12+ i\g) \big| / \log{\g} )}{\sqrt{\frac12\log\log \g}}  \in [a, b] \Bigg\},
$$
where the $\g$ are listed in increasing order. 
Our first theorem, a slight modification of  a recent result of \c{C}i\c{c}ek~\cite{log zeta}, provides the counting functions of these sequences. 

\begin{theorem}\label{thm 1}
Assume the Riemann hypothesis is true and that Hypothesis $\mathscr{H}_\d$ holds for some $\d\in(0, 1]$. Let $\max(|a|, |b|)\ll (\log \log \log T)^{\frac12-\e}$, where $\e>0$.
Then for all sufficiently large $T$,
\be\label{count formula 1}
\begin{split}
N_{[a, b]}(T) := \sum_{\substack{ 0<\g\leq T \\ \g\in  \G_{[a, b] }} }  1  =\frac{N(T)}{\sqrt{2\pi} }  \int_a^b e^{-{x^2}/{2}}\, dx +O\bigg( N(T)\frac{(\log\log\log T)^2}{\sqrt{\log\log T}}   \bigg) .    
\end{split}
\ee
For the sequence $\G_{[a, b]}^*$ we have
\be\label{count formula 2}
\begin{split}
N_{[a, b]}^*(T) := \sum_{\substack{ 0<\g\leq T \\ \g\in  \G_{[a, b]}^* } }  1 = 
 \frac{N(T) }{\sqrt{2\pi} }  \int_a^b e^{-{x^2}/{2}} dx  +o(N(T)).    
\end{split}
\ee
\end{theorem}
Observe that \eqref{count formula 2}  follows immediately from \eqref{count formula 1} since Hypothesis $\mathscr{H}_{\delta}$ implies that all but $o(N(T))$ of the zeros are simple.

 Our next theorem  is our main uniform distribution result.
\begin{theorem}\label{thm 2}
Assume the Riemann hypothesis and that Hypothesis $\mathscr{H}_\d$ is true for some $0<\d\leq 1$. 
Let    $a$ and $b$ be either fixed, or functions of $T$  for which $\max(|a|, |b|)\ll (\log \log \log T)^{\frac12-\e}$, where $\e>0$, and $\int_a^b  e^{-x^2/2} dx \gg 1$.
Then the sequences 
$\G_{[a, b] }$ and $\G^*_{[a, b] }$ are uniformly distributed modulo one.
\end{theorem}

The average gap between the ordinates   $ \g \in (0, T]$ is $2\pi/ \log T$ by  \eqref{zero count}. Thus  the numbers   $\g(\log T)/2\pi$ have average spacing  one. Not   surprisingly, it is more difficult to prove that these numbers are equidistributed modulo one. In fact, it  is not known.  This is also true for the numbers $\g(\log T)/2\pi$ with 
$ \g \in (0, T]$ and $\g\in\G_{[a,b]}$ or $\G_{[a,b]}^*$.  However,  if we  assume the following further conjecture  of the second author (confer \cite{GonekLandaulemma2}),  we can show uniform distribution  in all three cases.

\begin{conj}\label{Conjecture}
For $x, T \geq 2$ and any fixed $\epsilon >0$,
\begin{align*}
\sum_{0<\gamma \leq T} x^{i\gamma}  
&\ll \; T x^{-\frac12+\epsilon}+  T^{\frac12} x^\epsilon.
\end{align*}
\end{conj}

 Theorem 4 of  \cite{GonekLandaulemma2} provides evidence for Conjecture~\ref{Conjecture}.  It says that if $\psi(y)=\sum_{n\leq y} \Lambda(n)$,  where $\Lambda(n)$ is the von Mangoldt function, then Conjecture~\ref{Conjecture} implies that 
\be\label{prms short}
\psi(y+h) -\psi(y) =h +O(h^{\frac12} y^\epsilon)
\ee
for $1\leq h\leq y$ and $\epsilon>0$. Conversely, \eqref{prms short} implies
a weighted form of Conjecture~\ref{Conjecture}, namely,
\[
\sum_{\g} x^{i\g} \Big( \frac{\sin \g/2T}{\g/2T }\Big)^2
\ll \; T x^{-\frac12+\epsilon}+   T^{\frac12}  x^\epsilon.
\]

Using Conjecture~\ref{Conjecture}, one may  prove the following two theorems. In both, $\{x\}$ denotes the fractional part of $x$.
 \begin{theorem} \label{thm 3}
Assume the Riemann hypothesis and Conjecture~\ref{Conjecture}. If  $[\alpha,\beta]$ is a subinterval of $[0,1]$, then
 \be\label{discrp 0}
\sup_{\a, \b}\ \Bigg|  \; \sum_{ \substack{0<\g \leq T  \\ \{ \g (\log T)/2\pi\} \in [\alpha,\beta]}} 1 
\ - \  (\beta-\alpha)  N(T)\Bigg| 
=o(N (T)).
\ee
\end{theorem}

\begin{theorem} \label{thm 4}
Assume the Riemann hypothesis, Hypothesis $\mathscr{H}_\d$  for some $0<\d\leq 1$, and Conjecture~\ref{Conjecture}. 
Let  $a$ and $b$ be either fixed, or functions of $T$  for which $\max(|a|, |b|)\ll (\log \log \log T)^{\frac12-\e}$, where $\e>0$, and $\int_a^b  e^{-x^2/2} dx \gg 1$.  Then if $[\alpha,\beta]$ is a subinterval of $[0,1]$,  
\begin{equation}\label{discrp 1}
\sup_{\a, \b}\ \Bigg|  \; \sum_{\substack{0<\g \leq T, \,\g\in \G_{[a, b]} \\ \{ \g (\log T)/2\pi\} \in [\alpha,\beta]}} 1 
\ - \  (\beta-\alpha)  N_{[a, b]}(T)\Bigg| 
=o(N_{[a, b]}(T)).
\end{equation}
This also holds with $ \G_{[a, b]}$  replaced by $ \G_{[a, b]}^*$.
\end{theorem}

The method we use to prove Theorems~\ref{thm 1},  ~\ref{thm 2}, and  ~\ref{thm 4}  builds on techniques used in the first author's recent thesis to prove 
a discrete  analogue of Selberg's central limit theorem  (see~\cite{log zeta} and Lemma~\ref{lem 1} below).

Throughout we write $e(u)= e^{2\pi i u}$. We let $C$ denote a positive constant that may be different at different occurrences, and we  let 
$\mathbbm{1}_{[a,b]} $ denote the indicator function of the interval $[a, b]$.

\section*{Acknowledgements}
We thank the anonymous referee for their careful reading of our paper and for pointing out an error in the original version.
The first author thanks the Indian Institute of Technology Gandhinagar, where this work was begun, for its hospitality.
She is also grateful for a postdoctoral fellowship from the Pacific Institute for the Mathematical Sciences at the University of Northern British Columbia, where this work was completed.


\section{Proof of Theorem \ref{thm 1}}

Theorem~\ref{thm 1} follows easily from the next lemma.
\begin{lem}\label{lem 1}
Assume the Riemann hypothesis and that Hypothesis $\mathscr{H}_\d$ holds for some $\d\in(0, 1]$. Let
$m_\g$ denote the multiplicity of the zero $\r=\frac12+i\g$.
Then for all sufficiently large $T$,
\be\label{count lemma}
\begin{split}
 \# \bigg\{  0<\g\leq T \, :\,  & \frac{ \log(\big| \zeta^{(m_\g)} (\frac12+ i\g) \big| / (\log{T} )^{m_\g})}{\sqrt{\frac12\log\log T}} \in [a, b]  \bigg\} \\
=& \frac{ N(T)}{\sqrt{2\pi} }  \int_a^b e^{-{x^2}/{2}} \, dx +O\bigg(  N(T)\frac{(\log\log\log T)^2}{\sqrt{\log\log T}}   \bigg) .    
\end{split}
\ee
\end{lem}
This follows from the proof of Theorem 1.4 combined with Corollary 2.3 of  \c{C}i\c{c}ek~\cite{log zeta}.
 
Note  that
\be\notag
N_{[a, b]}(T)  
= \sum_{0 < \g \leq T}  \mathbbm{1}_{[a,b]} 
\bigg(\frac{\log \big(\big|\zeta^{(m_\g)}(\frac12+i\g) \big|/(\log \g)^{m_\g} \big) }{\sqrt{\frac12\log\log \g}}\bigg). 
\ee
Thus, to prove \eqref{count formula 1}, we  need to show that we may replace $\log \g$ and $\log\log \g$ here 
by $\log T$ and $\log\log T$, respectively, at the cost of a reasonable error term. To see this, note that
by \eqref{zero count}, the terms in the sum with $0 <\g\leq  T/\log T$ contribute at most  $O(T)$, hence
 \be\label{set S count}
N_{[a, b]}(T) =  \sum_{T/\log T < \g \leq T}  \mathbbm{1}_{[a,b]} 
\Bigg(\frac{\log \big(\big|\zeta^{(m_\g)}(\frac12+i\g) \big|/(\log \g)^{m_\g} \big) }{\sqrt{\frac12\log\log \g}}\Bigg)
+O(T). 
\ee
 For $T/\log T < \g \leq T$ we   easily find that
 \be\notag
\begin{split}
 \frac{\log \big(\big|\zeta^{(m_\g)}(\frac12+i\g) \big|/(\log \g)^{m_\g}  \big)  }{\sqrt{\frac12\log\log \g} }    
=  &\frac{ \log \big(\big|\zeta^{(m_\g)}(\frac12+i\g) \big|/(\log T)^{m_\g}  \big)  +O(m_\g \log\log T/ \log T)  }{\sqrt{\frac12\log\log T} }\\
   &\hskip1in \cdot\bigg(1 + O\Big(\frac{1}{\log T }\Big) \bigg).
\end{split}
\ee 
Using  \eqref{zero count} again, we see that $m_\g\ll \log T/\log\log T$. Thus,  if we impose the condition that
$\max(|a|, |b|)\ll (\log \log \log T)^{\frac12-\e}$, then when the expression on the left lies in the interval $[a, b]$, the right-hand side equals
 $$
 \frac{ \log \big(\big|\zeta^{(m_\g)}(\frac12+i\g) \big|/(\log T)^{m_\g}  \big)  }{\sqrt{\frac12\log\log T} }
+ O\Big(\frac{1}{ \sqrt{\log\log T} } \Big)  .
 $$
Using this with \eqref{count lemma}, we see that  replacing $\log \g$ and $\log\log \g$ in \eqref{set S count} by $\log T$ and $\log\log T$, respectively, changes \eqref{set S count} by  no more than
$\displaystyle
O( N(T) {(\log\log\log T)^2}/{\sqrt{\log \log T} }).
$ 
Hence
 \be\notag
N_{[a, b]}(T) =  \sum_{T/\log T < \g \leq T}  \mathbbm{1}_{[a,b]} 
\bigg(\frac{\log \big(\big|\zeta^{(m_\g)}(\frac12+i\g) \big|/(\log T)^{m_\g} \big) }{\sqrt{\frac12\log\log T}}\bigg)
+O\Big(N(T)\frac{(\log\log\log T)^2}{\sqrt{\log \log T} }\Big). 
\ee
Extending the sum back to the full range $0<\g\leq T$ and using  \eqref{count lemma} again, we see that
for sufficiently large $T$,
\be\notag\label{N_a,b}
N_{[a, b]}(T) =\frac{N(T)}{\sqrt{2\pi} }  \int_a^b e^{-{x^2}/{2}}\, dx  +O\Big( N(T)\frac{(\log\log\log T)^2}{\sqrt{\log \log T} }\Big),
\ee
provided $\max(|a|, |b|)\ll (\log \log \log T)^{\frac12-\e}$. This proves \eqref{count formula 1}. It has already been noted that 
\eqref{count formula 2} follows from \eqref{count formula 1} and Hypothesis $\mathscr{H}_{\delta}$, so the proof of Theorem~\ref{thm 1} 
is complete.

\section{Proof of Theorem~\ref{thm 2}}

We assume the Riemann hypothesis, Hypothesis $\mathscr{H}_\d$, and that $\max(|a|, |b|)\ll \\ (\log \log \log T)^{\frac12-\e}$ with $\e>0$. Our assumption that
$a$ and $b$ are either fixed, or functions of $T$  for which  $\int_a^b  e^{-x^2/2} dx \gg 1$ means,
 by Theorem~\ref{thm 1}, that $N_{[a, b]}(T) \gg N(T)$.
 Hence, by Weyl's criterion~\cite{Weyl},   the sequence   $ \G_{[a, b]}$ is uniformly distributed modulo one if, 
 for each fixed positive integer $\ell$,
 \be\label{Weyl 1} 
 \sum_{0 < \g \leq T} e(\ell \g) \,\, \mathbbm{1}_{[a,b]} 
\Bigg(\frac{\log \big(\big|\zeta^{(m_\g)}(\frac12+i\g) \big|/(\log \g)^{m_\g} \big) }{\sqrt{\frac12\log\log \g}}\Bigg) 
=o(N(T))
\ee
as $T\to\infty$. 
By the same argument we  used in the last section, replacing $\log \g$ and $\log\log \g$  here  by $\log T$ and $\log\log T$, respectively, changes the sum by at most $o(N(T)$. 
Thus, it suffices to show that 
 \be\notag\label{Weyl 2}
 \begin{split}
  \sum_{0 < \g \leq T} e(\ell \g) \,\, \mathbbm{1}_{[a,b]} 
\Bigg(  \frac{\log \big(\big|\zeta^{(m_\g)}(\frac12+i\g) \big|/(\log T)^{m_\g} \big) }{\sqrt{\frac12\log\log T}}\Bigg)  
 = o( N(T)).
\end{split}
\ee
  
Now let
$ \displaystyle P(\g) =\sum_{p \leq X^2} \frac{1}{ p^{1/2+i\g} },$
where $p$  runs over the primes. 
 The gist of  Corollary 2.3 in \cite{log zeta} and some of the analysis following it,
 is that $P(\g)$ is  on average a good approximation   to 
 $$
 \frac{\log \big(\big|\zeta^{(m_\g)}(\frac12+i\g) \big|/(\log T)^{m_\g} \big) }{\sqrt{\frac12\log\log T}},
 $$
 provided $X$ is sufficiently large.
Indeed, from the discussion in Sections 5.5 and 6 of \cite{log zeta} it follows that
  \begin{align*}
\sum_{0<\g \leq T}\mathbbm{1}_{[a,b]} 
\Bigg( \frac{\log \big(\big|\zeta^{(m_\g)}(\frac12+i\g) \big|/(\log T)^{m_\g} \big) }{\sqrt{\frac12\log\log T}}  \Bigg)
= \sum_{0<\g \leq T} & \mathbbm{1}_{[a,b]} \Bigg(\frac{\Re {P(\g)}}{\sqrt{\frac12\log\log T}}\Bigg) \\
         &+O  \bigg(N(T) \frac{(\log\log\log T)^2}{\sqrt{\log\log T}}   \bigg) .
 \end{align*} 
This is the key result that allows us to prove our theorem. An immediate consequence is that 
\[
\mathbbm{1}_{[a,b]} \Bigg(\frac{\log \big(\big|\zeta^{(m_\g)}(\frac12+i\g) \big|/(\log T)^{m_\g} \big) }{\sqrt{\frac12\log\log T}}\Bigg)
= \mathbbm{1}_{[a,b]} \Bigg(\frac{\Re{P(\g)}}{\sqrt{\frac12\log\log T}}\Bigg)
\]
for all but $\displaystyle O\big(N(T)  (\log\log\log T)^2/\sqrt{\log\log T}   \big)=o(N(T))$ values  of $\g$ in $(0, T]$.  
Therefore, 
\be
\begin{split}\label{Weyl 3}
 \sum_{0<\g \leq T} & e(\ell\g) \mathbbm{1}_{[a,b]} \Bigg(\frac{\log \big(\big|\zeta^{(m_\g)}(\frac12+i\g) \big|/(\log T)^{m_\g} \big) }{\sqrt{\frac12\log\log T}}\Bigg)  
 \\
=& \sum_{0<\g \leq T} e(\ell\g) \mathbbm{1}_{[a,b]}  \Bigg(\frac{\Re P(\g)}{\sqrt{\frac12\log\log T}}\Bigg)
+ O  \bigg(N(T) \frac{(\log\log\log T)^2}{\sqrt{\log\log T}}   \bigg).
\end{split}
\ee

Writing
\be\notag
A=A(T)=a\sqrt{\frac12\log\log T} \quad \text{and} \quad 
B=B(T)=b\sqrt{\frac12\log\log T},
\ee
we see that to prove our theorem we must show that \be\label{exp ell sum}
\sum_{0<\g \leq T} 
e(\ell\g)\mathbbm{1}_{[A, B]}  ( \Re P(\g)) = o(N(T))
\ee
for each  positive integer $\ell$.
To do this, we  replace the characteristic function $\mathbbm{1}_{[A, B]}$ by an approximation. 
Let $\Omega > 0$ and set
    \be\label{eq:F}
    F_{\Omega}(x)
    =\Im \int_0^\Omega G\Big(\frac \omega \Omega\Big)\exp{(2\pi ix\omega)}\frac{\mathop{d\omega}}{\omega},
    \ee
where
\be\notag
    G(u)= \frac{2u}{\pi}+2u(1-u)\cot{(\pi u)} \quad \text{for} \quad u\in[0, 1]. 
\ee
Then 
    \be\label{eq:sgn}
    \sgn{(x)}
    =F_{\Omega}(x)
    +O\bigg(\sdfrac{\sin^2(\pi\Omega x)}{(\pi\Omega x)^2}\bigg).
    \ee
(see~\cite[pp. 26--29]{Tsang}). It follows that 
    \be\label{eq:chi F}
    \mathbbm{1}_{[A, B]}(x)
    =\frac12F_\Omega(x-A)-\frac12F_\Omega(x-B)
    +O\bigg(\sdfrac{\sin^2(\pi\Omega(x-A))}{(\pi\Omega(x-A))^2}\bigg)
    +O\bigg(\sdfrac{\sin^2(\pi\Omega(x-B))}{(\pi\Omega(x-B))^2}\bigg).
    \ee
    This is the desired approximation of $ \mathbbm{1}_{[A, B]}$.
Here we take $\displaystyle x=\Re{ {P}(\g)} = \Re \sum_{p \leq X^2} \frac{1}{ p^{1/2+i\g} }$ and 
 \be\label{X, Omega}
 X = T^{\frac{1}{(\log\log T)^{20}}}, \qquad   \Omega=(\log\log T)^2.
 \ee  
  Now, it was shown in the course of the proof of Proposition 5.5 in~\cite{log zeta} (with slightly different notation and parameters) that
\be\notag
\sum_{0 < \g \leq T}\frac{\sin^2\big(\pi \Omega(\Re {P}(\g)-A)\big)}{\big(\pi \Omega(\Re {P}(\g)-A)\big)^2}
 \ll  \frac{N(T)}{\Omega},  
\ee
and similarly for the sum with $A$ replaced by $B$.  
Thus, by \eqref{eq:chi F},
\be\label{set up}
\begin{split}
& \sum_{0<\g \leq T}  
e(\ell\g) \mathbbm{1}_{[A, B]}  ( \Re {P}(\g)) \\
 = &\, \sdfrac12  \sum_{0<\g \leq T}   e(\ell\g)  F_{\Omega}  (\Re{ P(\g)}-A)
-\sdfrac12 \sum_{0<\g \leq T}  e(\ell\g) F_{\Omega}(\Re {P}(\g)-B) 
+O\Big( \frac{N(T)}{\Omega}\Big) .
\end{split}
\ee
From this and \eqref{exp ell sum} we  see that it suffices to
prove that 
\be\label{sum F_A}
\sum_{0 < \g \leq T} e(\ell \g) F_{\Omega}(\Re  {P}(\g)-A) =o(N(T))
\ee
for each positive integer $\ell$, and similarly for the sum with   $A$ replaced by $B$.

To this end we use  \eqref{eq:F} to write
    \be\label{eq:substitute F}
    \begin{split}
    \sum_{0 < \g \leq T} e(\ell \g)   F_{\Omega}(\Re {P}(\g)-A)    
    =&\sum_{0 < \g \leq T}  e(\ell \g)  \Im \int_0^\Omega G\Big(\frac \omega \Omega\Big)
    e^{-2\pi i A \omega} \exp{\big(2\pi i\omega\Re  {P}(\g) \big)}\frac{\mathop{d\omega}}{\omega}.
    \end{split}
    \ee
By Taylor's theorem, for any positive integer $K$,
    \be\notag
    \exp\bigs(2\pi i\omega\Re{P}(\g)\bigs)
    =1+\sum_{1\leq k < K} \frac{(2\pi i\omega\Re {P}(\g))^k}{k!} + 
    O\Big(\frac{(2\pi \omega|\Re {P}(\g)|)^K}{K!}\Big).
    \ee
Inserting this  in \eqref{eq:substitute F} and taking 
\be\label{K}
K=2\big[(\log\log T)^6\big],   
\ee 
where $[ x]$ denotes the greatest integer less than or equal to $x$,
we obtain
      \begin{align} \notag
    \sum_{0 < \g \leq T} e(\ell \g)  F_{\Omega}(\Re {P}(\g)-A)
    =&\, F_\Omega(A)  \sum_{0 < \g \leq T} e(\ell \g) \\
    \label{exp}%
    +\sum_{0 < \g \leq T} e(\ell \g)  \Im \int_0^\Omega & G\Big(\frac \omega \Omega\Big) 
    e^{-2\pi i A \omega} \sum_{1\leq k < K}\frac{(2\pi i\omega)^k}{k!} 
     (\Re {P}(\g))^k  \frac{\mathop{d\omega}}{\omega}\\
   \notag
    &+ O\bigg( \sum_{0 < \g \leq T} |\Re {P}(\g)|^K \int_0^\Omega  G\Big(\frac \omega \Omega\Big) \frac{(2\pi\omega)^K}{K!}
    \frac{\mathop{d\omega}}{\omega} \bigg) .
    \end{align}
We estimate  the  sums over $\g$ on the right-hand side of the equation by means of the following result, which is an immediate consequence of an
unconditional theorem and its corollary in~\cite{GonekLandaulemma2} (see~\cite{GonekLandaulemma1} also). 
\begin{lem}\label{Land-Gon} Assume the Riemann hypothesis and  let $x, T>1$.
 Then
 \be\label{L-G}
 \begin{split}
 \sum_{0<\g\leq T} x^{i\g}= - &\frac{T}{2\pi}\frac{\Lambda(x)}{\sqrt x} +O(\sqrt{x} \log 2xT \log\log 3x) 
 +O\Big(\log x\min\Big(\frac{T}{\sqrt x}, \frac{\sqrt x}{\langle x \rangle}\Big)\Big) \\
 &+O\Big(\log 2T \min\Big(\frac{T}{\sqrt x}, \frac{1}{\sqrt x \log x} \Big)\Big).
 \end{split}
 \ee
Here $\Lambda(x)=\log p$ if $x$ is a positive integral power of a prime $p$ and $\Lambda(x)=0$
for all other real numbers $x$, and
$\langle x\rangle$ denotes the distance from $x$ to the nearest prime power other than $x$ itself.
If $0<x<1$, \eqref{L-G} also holds provided we replace $x$ on the right-hand side by $1/x$.
\end{lem}
When $x>1$, we will  write \eqref{L-G}  as
\be\label{L-G short 1}
 \sum_{0<\g\leq T} x^{i\g} =M(x) + E_1(x) + E_2(x)+ E_3(x),
 \ee
where  $M(x)$ and  $E_i(x), \ i=1, 2, 3,$  also depend on $T$.  When  $0<x<1$,   all  the $x$'s  on the  
 right-hand side of  \eqref{L-G short 1} 
are to be replaced by $1/x$. Note that  when $x=1$,  \; $ 
 \sum_{0<\g\leq T} x^{i\g}=N(T)$.

Returning to  \eqref{exp},  observe that by \eqref{eq:sgn}, $F_{\Omega}(A) \ll 1$.
Furthermore, taking $x=e^{2\pi \ell }>1$ in \eqref{L-G},  we find that $\sum_{0<\g\leq T}e(\ell \g) \ll T$.
Thus, the first term on the right-hand side of \eqref{exp} is
\be\label{1st term}
F_\Omega(A) \sum_{0 < \g \leq T} e(\ell \g)  \ll T.
\ee

For the final term in \eqref{exp} we use  Lemma 5.2 of~\cite{log zeta}, which says that
\be\notag
\sum_{0<\g \leq T} | \Re P(\g) |^K \ll (cK\Psi)^{K/2} N(T),
\ee
where 
$
     \Psi =\log\log T.
$
From this and Stirling's approximation, we find that the $O$-term in \eqref{exp} is
    \begin{align*}
      &\ll N(T) \int_0^\Omega  G\Big(\frac \omega \Omega\Big)   \frac{(2\pi \omega)^K}{K!}(cK\Psi)^{K/2} \frac{\mathop{d\omega}}{\omega} \\
       &\ll N(T) \int_0^\Omega  G\Big(\frac \omega \Omega\Big)  \frac{\omega(2\pi e)^K\omega^{K-1}}{K^K} (cK\Psi)^{K/2}
       \frac{\mathop{d\omega}}{\omega} .
     \end{align*}
By \eqref{X, Omega} and \eqref{K}, and since $G$ is bounded, this is
    \[
       \ll N(T)\int_0^\Omega  G\Big(\frac \omega \Omega\Big)   \bigg(\frac{c\, \Omega\sqrt{\Psi}}{\sqrt{K}}\bigg)^K 
      \mathop{d\omega}\\
      \ll \frac{N(T) }{2^{K}}\int_0^\Omega  G\Big(\frac \omega \Omega\Big)   \mathop{d\omega}
     \ll T.
    \]
Combining this and \eqref{1st term}, we may rewrite \eqref{exp}  as
    \be\label{F Omega}
    \begin{split}
    & \sum_{0 < \g \leq T}  e(\ell \g)  F_{\Omega}(\Re P(\g)-A) \\
     =& \int_0^\Omega  G\Big(\frac \omega \Omega\Big) 
    \sum_{1\leq k < K}   \Im{\big(e^{-2\pi i A \omega}\, i^k\big)} \frac{(2\pi\omega)^k}{k!} 
    \sum_{0 < \g \leq T} e(\ell \g)  (\Re P(\g))^k
     \frac{\mathop{d\omega}}{\omega} 
    +O(T).
    \end{split}
    \ee    
To estimate the  right-hand side we next bound the sums
\be\label{S(k)}
  S(k) =  \sum_{0 < \g \leq T} e(\ell \g)  (\Re P(\g))^k.
\ee
By the binomial theorem
    \[
    S(k)  = \frac{1}{2^k}\sum_{j=0}^k \binom{k}{j} 
    \sum_{0 < \g \leq T} e(\ell \g)  \Big(\sum_{p\leq X^2}\frac{1}{p^{1/2+i\g}}\Big)^j
    \Big(\sum_{p\leq X^2}\frac{1}{p^{1/2-i\g}}\Big)^{k-j}.
    \]
Let $a_{r}(p_1\dots p_{r})$ denote the number of permutations of the primes $p_1,\dots, p_{r}$, which might or might not be distinct. Also, for the rest of the paper, $n$ will always denote a product of $j$ primes, each of which is at most $X^2$, while $m$ denotes a product of $k-j$ primes, again each of size at most $X^2$. We may thus write
    \begin{align*}
    S(k) =&\, \frac{1}{2^k}\sum_{j=0}^k \binom{k}{j} 
    \sum_{0 < \g \leq T} e(\ell \g)  \sum_n \frac{a_j(n)}{n^{1/2+i\g}} \sum_m \frac{a_{k-j}(m)}{m^{1/2-i\g}} \\
    =&\,\frac{1}{2^k}\sum_{j=0}^k \binom{k}{j}  \sum_n \frac{a_j(n)}{\sqrt n} \sum_m \frac{a_{k-j}(m)}{\sqrt m}
    \sum_{0 < \g \leq T} \Big(\frac{m e^{2\pi\ell}}{n}\Big)^{i\g}.
    \end{align*}
Since $e^{2\pi \ell} = (-1)^{-2i\ell}$ is of the form $\a_0^{\b_0}$ with $\a_0, \b_0$ algebraic, $\a_0\neq 0,1$ and $-2i\ell$ not rational, the Gelfond-Schneider theorem implies that $e^{2\pi \ell}$ is transcendental. Thus, $m e^{2\pi\ell}/{n}$  can neither  be a positive integer nor the reciprocal of a positive integer. The $M$ term in \eqref{L-G short 1} is therefore always zero. Hence, we may write
 \be\label{contribution for k}
    \begin{split}
S(k) &= \frac{1}{2^k}\sum_{j=0}^k \binom{k}{j}  \sum_n \frac{a_j(n)}{\sqrt n}
 \sum_{\substack{m\\{me^{2\pi\ell}}/{n} >1}}  \frac{a_{k-j}(m)}{\sqrt m} 
 \bigg( \sum_{i=1}^3 E_i \Big(\frac{m e^{2\pi\ell}}{n}\Big)\bigg)\\
 &\quad + \frac{1}{2^k}\sum_{j=0}^k \binom{k}{j}  \sum_n \frac{a_j(n)}{\sqrt n}
 \sum_{\substack{m\\{me^{2\pi\ell}}/{n} <1}}  \frac{a_{k-j}(m)}{\sqrt m} 
 \bigg( \sum_{i=1}^3 E_i \Big(\frac{n}{m e^{2\pi\ell}}\Big)\bigg)
 \\
& =: S_1(k) + S_2(k) . 
  \end{split}
    \ee
To estimate $S_1(k)$, we insert the bounds for $E_1, E_2,$ and $E_3$ from \eqref{L-G} in to obtain
\be\label{S1}
S_1(k) =\mathcal E_1(k)+\mathcal E_2(k)+\mathcal E_3(k), 
\ee
where
    \be\notag 
    \begin{split}
      \mathcal E_1(k)&= \frac{1}{2^k}\sum_{j=0}^k \binom{k}{j}  \sum_n \frac{a_j(n)}{\sqrt n} 
     \sum_{\substack{m\\{me^{2\pi\ell}}/{n} >1}}
     \frac{a_{k-j}(m)}{\sqrt m}  \sqrt{\frac{me^{2\pi \ell}}{n}} \log\Big(\frac{me^{2\pi \ell}T}{n}\Big) \log\log\Big(\frac{3me^{2\pi \ell}}{n}\Big),  \\
   \mathcal  E_2(k)&= \frac{1}{2^k}\sum_{j=0}^k \binom{k}{j}  \sum_n \frac{a_j(n)}{\sqrt n}  \sum_{\substack{m\\{me^{2\pi\ell}}/{n} >1}} \frac{a_{k-j}(m)}{\sqrt m}  \frac{ \log(me^{2\pi \ell}/n)}{\sqrt{me^{2\pi \ell}/n}}  
    \min\bigg(T, \frac{me^{2\pi \ell}/n }{\langle me^{2\pi \ell}/n  \rangle}  \bigg), \\
  \mathcal   E_3(k)&= \frac{1}{2^k}\sum_{j=0}^k \binom{k}{j}  \sum_n \frac{a_j(n)}{\sqrt n}  \sum_{\substack{m\\{me^{2\pi\ell}}/{n} >1}} \frac{a_{k-j}(m)}{\sqrt m} \frac{\log  T}{\sqrt{me^{2\pi \ell}/n} }  \min \bigg(T, \frac{ 1}{\log(me^{2\pi \ell}/n)}  \bigg).
    \end{split}
    \ee

First consider $\mathcal E_1(k)$. Since $e^{2\pi \ell}$ is fixed, we see that
    \begin{align*}
    \mathcal E_1(k)
    &= \frac{1}{2^k}\sum_{j=0}^k \binom{k}{j}  \sum_n \frac{a_j(n)}{n} \sum_{\substack{m\\{me^{2\pi\ell}}/{n} >1}} a_{k-j}(m) \log\Big(\frac{mT}{n}\Big) \log\log \Big(\frac{3m}{n}\Big) \\
    &\ll  \frac{\log T\log\log T }{2^k}\sum_{j=0}^k \binom{k}{j}  \sum_n  {a_j(n)} 
    \sum_{\substack{m\\{me^{2\pi\ell}}/{n} >1}} a_{k-j}(m). 
    \end{align*}
Here we have dropped the $n$ in the denominator in the sum over $n$ and used \eqref{X, Omega} and \eqref{K} to deduce that $m \leq X^{2k}  \leq X^{2K}  \leq T$. Next, from the 
definitions of $a_j(n)$ and $a_{k-j}(m)$ and by the binomial theorem, we see that
    \be \label{eq:psi pi}
    \begin{split}
    & \frac{1}{2^k}\sum_{j=0}^k \binom{k}{j}  \sum_n  a_j(n) 
      \sum_{\substack{m\\{me^{2\pi\ell}}/{n} >1}} a_{k-j}(m) 
      \\
    &\ll   \frac{1}{2^k}\sum_{j=0}^k \binom{k}{j} 
    \Big( \sum_{p\leq X^2}1\Big)^j  \Big( \sum_{p\leq X^2} 1\Big)^{k-j}  
 = \pi(X^2)^k,
    \end{split} 
    \ee
where $\pi(X^2)$ denotes the number of primes up to $X^2$. By the prime number theorem 
$\pi(X^2)\ll  {X^{2}}/{\log X}, $
so 
    \be\notag
\mathcal    E_1(k) \ll \log T \log\log T \frac{X^{2k}}{(\log X)^k}.
    \ee

To estimate   
\be\label{E2}
     \mathcal  E_2(k)= \frac{1}{2^k}\sum_{j=0}^k \binom{k}{j}  \sum_n \frac{a_j(n)}{\sqrt n}  \sum_{\substack{m\\{me^{2\pi\ell}}/{n} >1}} \frac{a_{k-j}(m)}{\sqrt m}  \frac{ \log(me^{2\pi \ell}/n)}{\sqrt{me^{2\pi \ell}/n}}  
    \min\bigg(T, \frac{me^{2\pi \ell}/n }{\langle me^{2\pi \ell}/n  \rangle}  \bigg),
\ee
we require a lower bound for $\displaystyle \langle  {me^{2\pi\ell}}/{n}\rangle$. Recall that $\langle x\rangle $ denotes  the distance from $x$ to the nearest prime power other than $x$ itself. As $e^{2\pi\ell}$ is transcendental, $me^{2\pi \ell}/n$ is not an integer. 
Now, for any positive non integral   real number $x$, we have
$$
 \langle  x \rangle \geq \min_{r\in \mathbb Z} |x-r| = \min\{ x-[x], 1-x+[x]  \}.
$$
Thus 
$$
\Big\langle \frac{me^{2\pi\ell}}{n}\Big\rangle 
\geq \min\bigg\{     \frac{me^{2\pi\ell}}{n}- \Big[\frac{me^{2\pi\ell}}{n} \Big]  ,   \;
1-\frac{me^{2\pi\ell}}{n}  +  \Big[\frac{me^{2\pi\ell}}{n} \Big]  \bigg\}.
$$
We now recall a special case of Baker's theorem~\cite[p. 24]{Baker}. Since 
$e^{2\pi \ell} = (-1)^{-2i \ell}$,  for a given positive integer $\ell$,
    \be\notag
    \Big|e^{2\pi \ell}-\frac pq \Big| >q^{-C \log\log q}
    \ee
    for all rationals $p/q$ ($p\geq 0, q\geq 4$),
where $C$ is a constant that  depends on $\ell$. 
Thus  
\be\label{Baker}
    \Big| \frac{me^{2\pi \ell}}{n} - \frac mn\frac pq \Big| > \frac mn q^{-C \log\log q}.
    \ee
If we let $q=m$ and $\displaystyle  p=n\Big[\frac{me^{2\pi\ell}}{n} \Big] $, we obtain
$$
  \frac{me^{2\pi\ell}}{n}- \Big[\frac{me^{2\pi\ell}}{n} \Big] 
 \geq \frac mn \, m^{-C \log\log m}.
$$
Similarly, taking $q=m$ and $\displaystyle  p=n\Big[\frac{me^{2\pi\ell}}{n} \Big] +n$, we see that
$\displaystyle
\Big( 1-  \frac{me^{2\pi\ell}}{n} + \Big[\frac{me^{2\pi\ell}}{n} \Big] \Big)
$
has the same lower bound. Hence
\be\label{prime distance}
\Big\langle \frac{me^{2\pi\ell}}{n}\Big\rangle 
 \geq \frac mn \, m^{-C \log\log m},
\ee
provided $m\geq 4$.
Now the terms in \eqref{E2} with $1\leq m \leq 3$ contribute  
 \be\notag
\begin{split}
 &   \ll_{\ell} 
     \frac{1}{2^k}\sum_{j=0}^k \binom{k}{j}  \sum_{n<6 e^{2\pi\ell}} \frac{a_j(n)}{n} \sum_{m\leq 3} a_{k-j}(m)
   \\
    &\ll \frac{1}{2^k}\sum_{j=0}^k \binom{k}{j} \Big(\sum_{p\leq  6 e^{2\pi\ell} }1 \Big)^j \Big(\sum_{p\leq 3}1 \Big)^{k-j} 
     \ll e^{ C k},
\end{split}
\ee
where $C$ is a constant depending on $\ell$.
Using \eqref{prime distance} to estimate the terms in \eqref{E2} with $m\geq 4$, we find that they contribute  
 \be\notag
\begin{split}
&  \ll 
    \frac{\log T}{2^k}\sum_{j=0}^k \binom{k}{j}  \sum_n a_j(n) \sum_m  {a_{k-j}(m)} 
    m^{C \log\log m} \\
   & \ll\frac{\log T}{2^k}  X^{3Ck \log\log X}\sum_{j=0}^k \binom{k}{j} \Big(\sum_{p\leq X^2 }1\Big)^j \Big(\sum_{p\leq X^2 }1\Big)^{k-j} \\
    & \ll  ( \log T ) \pi(X^2 )^{k}  X^{3Ck \log\log X}
     \ll (\log T )  X^{4Ck \log\log X}.
\end{split}
\ee
Thus, 
\be\notag
\mathcal E_2(k) \ll_{\ell} (\log T )  X^{4Ck \log\log X} +e^{C k}  \ll_{\ell} (\log T )  X^{4Ck \log\log X}.
\ee

Next consider  
\be \label{E 3}
   \mathcal E_3(k) =\frac{\log T}{2^k}\sum_{j=0}^k \binom{k}{j}  \sum_n  a_j(n)  \sum_{\substack{m\\{me^{2\pi\ell}}/{n} >1}} \frac{a_{k-j}(m)}{m}  \min \bigg(T, \frac{ 1}{\log(me^{2\pi\ell }/n)}  \bigg). 
\ee
 Since $m, n \leq X^{2k}$, by \eqref{Baker} with $p=n, q=m$ and $m\geq 4$, we have
 \be\label{Baker 2}
 \frac{me^{2\pi \ell}}{n} - 1  > \frac mn m^{-C \log\log m} \geq X^{-2k-2kC\log\log (X^{2k})} 
    \gg X^{-3kC\log\log X }.
    \ee    
Thus,     
$$
\frac{1}{\log(me^{2\pi\ell }/n)}\ll X^{3kC\log\log X}
$$   
for $m\geq 4$.    The terms in \eqref{E 3} with $m\leq 3$ contribute
\be\notag
\begin{split}
&\ll \frac{\log T}{2^k}\sum_{j=0}^k \binom{k}{j}  \sum_{n<6 e^{2\pi\ell} }  a_j(n)  \sum_{m\leq 3} {a_{k-j}(m)}  \\
&=\frac{\log T}{2^k}\sum_{j=0}^k \binom{k}{j}  \Big(\sum_{p\leq  6 e^{2\pi\ell} }1 \Big)^j \Big(\sum_{p\leq 3}1 \Big)^{k-j} 
\ll  e^{Ck} \log T.
\end{split}
\ee  
The contribution of the terms with $m\geq 4$ is
    \[
      \begin{split}
    \mathcal E_3(k)  
    &\ll  X^{3kC\log\log X} \ \frac{\log T }{2^k}  
    \sum_{j=0}^k \binom{k}{j}   \sum_n   a_j(n) \sum_{m}  a_{k-j}(m) \\
    &=X^{3kC\log\log X} \  \pi(X^2)^k \log T 
    \ll \frac{X^{2k} \log T}{(\log X)^{k}}  X^{3kC\log\log X}
    \ll X^{4kC\log\log X} \log T.
       \end{split}
    \]
 Thus
    \be\notag
  \mathcal  E_3(k) \ll  X^{4kC\log\log X} \log  T.
    \ee


Combining our estimates for  $\mathcal E_1(k), \mathcal E_2(k) $,  and   $\mathcal E_3(k)$  in \eqref{S1}, we see that 
$$ S_1(k) \ll_\ell  (\log T )  X^{4Ck \log\log X}.
$$
The estimation of $S_2(k)$ is very similar and leads to the same bound. Hence, by \eqref{contribution for k}
$$
S(k) \ll_\ell (\log T )  X^{4Ck \log\log X}.
$$
Inserting this into \eqref{F Omega},
 we find that
    \be\notag
    \begin{split}
     \sum_{0 < \g \leq T}  e(\ell \g)  F_{\Omega}(\Re P(\g)-A) 
    \ll_\ell   X^{4CK\log\log X} \log T\int_0^\Omega  G\Big(\frac \omega \Omega\Big) 
    \sum_{1\leq k < K}   \frac{(2\pi\omega )^k}{k!} 
     \frac{\mathop{d\omega}}{\omega} 
   +T.
     \end{split}
    \ee    
As $G$ is bounded over $[0, 1]$, this is
$$
\ll_\ell  \Omega \,e^{2\pi \Omega}  X^{4CK\log\log X} \log T  +T.
$$
By the choice of parameters $X, \Omega,$ and $K$ in  \eqref{X, Omega} and \eqref{K}, it follows that  for a fixed $\ell$,
$$
  \sum_{0 < \g \leq T}  e(\ell \g)  F_{\Omega}(\Re P(\g)-A) \ll e^{7(\log\log T)^2} T^{8C /(\log\log T)^{13}} +T 
  \ll  T =o(N(T)).
  $$
This establishes \eqref{sum F_A}, and 
the same estimate clearly holds when $A$ is replaced by $B$. 
This completes the proof   that 
$\G_{[a, b]}$ is uniformly distributed modulo one.

Since the number of $\g$'s in $\G_{[a, b]}$ with $0<\g\leq T$ that are not elements of $\G_{[a, b]}^*$ (in other words, that are not simple)
is at most $o(N(T))$, we see that $\G_{[a, b]}^*$ is also uniformly distributed modulo one.

 \section{Proof of Theorem~\ref{thm 3}}

By the Erd\H{o}s-Tur\'an inequality (see \cite{Mont}, Chapter 1, Corollary 1.1), if $L$ is a positive integer and $[\alpha,\beta]$ is a subinterval of $[0,1]$, then
\begin{equation}\label{ErdosTuran 1}
\Bigg| \sum_{\substack{0<\g \leq T, \\ \{ \g (\log T)/2\pi \} \in [\alpha,\beta]}} 1 \ \ - \ (\beta-\alpha)N(T) \Bigg| \ \leq \ \frac{N(T) }{L+1} \ + \ 3\sum_{\ell\leq L} \frac{1}{\ell}\Bigg| \sum_{0<\g\leq T } e\Big(\ell \g \frac{\log T}{2\pi}\Big)\Bigg|.
\end{equation} 
By Conjecture~\ref{Conjecture}, for each integer $\ell >0$,
\be\notag
\sum_{0<\g\leq T } e\Big(\ell \g \frac{\log T}{2\pi}\Big)= \sum_{0 < \g \leq T} T^{i\ell \gamma} \ll T^{\frac12+\epsilon \ell}.
\ee 
Hence, the right-hand side of \eqref{ErdosTuran 1} is 
$$
\ll \frac{N(T)}{L} + (\log L )\ T^{\frac12+\epsilon L}. 
$$
Taking $ \displaystyle L=\Big[\frac{1}{2\epsilon}\Big]$ and assuming that $0<\epsilon <1/2$, we see that this is 
$\ll \epsilon N(T)$. Since $\epsilon$ can be arbitrarily small, this
establishes \eqref{discrp 0}.

\section{Proof of Theorem~\ref{thm 4}}

By the Erd\H{o}s-Tur\'an inequality again, 
 if $L$ is a positive integer and $[\alpha,\beta]$ is a subinterval of $[0,1]$, then
\begin{equation}\label{ErdosTuran}
\begin{split}
\Bigg| \sum_{\substack{0<\g \leq T, \g\in \G_{[a, b]}\\ \{ \g (\log T)/2\pi \} \in [\alpha,\beta]}} 1 \ \ &- \ (\beta-\alpha)N_{[a, b]}(T) \Bigg| \ \\
\leq \ 
&\frac{N_{[a, b]}(T) }{L+1} \ + \ 3\sum_{\ell\leq L} \frac{1}{\ell}\Bigg| \sum_{0<\g\leq T, \g\in \G_{[a, b]}} e\Big(\ell \g \frac{\log T}{2\pi}\Big)\Bigg|.
\end{split}
\end{equation} 
 Thus, to prove \eqref{discrp 1},  we need to estimate
 \be\label{Weyl sum 1}
 \begin{split}
   \sum_{0<\g\leq T, \g\in \G_{[a, b]}} e\Big(\ell \g \frac{\log T}{2\pi}\Big)=& 
    \sum_{0 < \g \leq T} T^{i\ell  \g } \,\, \mathbbm{1}_{[a,b]} 
\Bigg(\frac{\log \big(\big|\zeta^{(m_\g)}(\frac12+i\g) \big|/(\log \g)^{m_\g} \big) }{\sqrt{\frac12\log\log \g}}\Bigg)
\end{split}
 \ee
for positive integers $\ell$. 
We do this, for the most part, by following the procedure of estimating the corresponding sum in \eqref{Weyl 1} in the previous section.
To start with, the same analysis that led to \eqref{Weyl 3} leads to 
 \be\label{Weyl sum 2}
 \begin{split} 
 &   \sum_{0 < \g \leq T} T^{i\ell  \g } \,\, \mathbbm{1}_{[a,b]} 
\Bigg(\frac{\log \big(\big|\zeta^{(m_\g)}(\frac12+i\g) \big|/(\log \g)^{m_\g} \big) }{\sqrt{\frac12\log\log \g}}\Bigg)
\\
=&\sum_{0<\g \leq T} T^{i\ell  \g } \,\, \mathbbm{1}_{[a,b]}  \Bigg(\frac{\Re P(\g)}{\sqrt{\frac12\log\log T}}\Bigg)
+ O \bigg(N(T) \frac{(\log\log\log T)^2}{\sqrt{\log\log T}}   \bigg).
\end{split}
 \ee
Similarly, the analysis that led to \eqref{set up}, with the same choices of the parameters 
$A, B, X$, and $\Omega$, and Dirichlet polynomial $P$, shows  that 
\be\label{set up 2}
\begin{split}
\sum_{0<\g \leq T}  
&T^{i\ell  \g } \mathbbm{1}_{[A, B]}  ( \Re {P}(\g)) \\
 = &\sdfrac12  \sum_{0<\g \leq T}  T^{i\ell  \g }   F_{\Omega}  (\Re{ P(\g)}-A)
-\sdfrac12 \sum_{0<\g \leq T}  T^{i\ell  \g }  F_{\Omega}(\Re {P}(\g)-B) 
+O\Big( \frac{N(T)}{\Omega}\Big) .
\end{split}
\ee
And, similarly to \eqref{exp}, we find that
  \begin{align} \notag
    \sum_{0 < \g \leq T} T^{i\ell  \g }   F_{\Omega}(\Re {P}(\g)-A)
    = &\, F_\Omega(A)  \sum_{0 < \g \leq T} T^{i\ell  \g }  \\
    \label{exp 5}
    +\sum_{0 < \g \leq T} T^{i\ell  \g }  \Im \int_0^\Omega & G\Big(\frac \omega \Omega\Big) 
    e^{-2\pi i A \omega} \sum_{1\leq k < K}\frac{(2\pi i\omega)^k}{k!} 
     (\Re {P}(\g))^k  \frac{\mathop{d\omega}}{\omega}\\
   \notag &+ O\bigg( \sum_{0 < \g \leq T} |\Re {P}(\g)|^K \int_0^\Omega  G\Big(\frac \omega \Omega\Big) \frac{(2\pi\omega)^K}{K!}
    \frac{\mathop{d\omega}}{\omega} \bigg) ,
    \end{align}
where $K=2\big[(\log\log T)^6\big]$.

By \eqref{eq:sgn}, $F_{\Omega}(A) \ll 1$ and, by Conjecture~\ref{Conjecture},
\be\notag
 \sum_{0 < \g \leq T} T^{i\ell \gamma} \ll T^{\frac12+\epsilon \ell}
\ee
for any $\epsilon>0$. Thus, the first term on the right-hand side of \eqref{exp 5} is $O(T^{\frac12+\epsilon \ell})$.
 The third term is estimated in the same way as the third term in \eqref{exp}  and is likewise $O(T)$.
Hence
    \be\label{F Omega 3}
    \begin{split}
  & \sum_{0 < \g \leq T}  T^{i\ell\gamma} F_{\Omega}(\Re{P}(\g)-A) 
   \\
    = &\int_0^\Omega   G\Big(\frac \omega \Omega\Big) 
    \sum_{1\leq k < K}   \Im{\big(e^{-2\pi i A \omega}\, i^k\big)} \frac{(2\pi\omega)^k}{k!} 
    \sum_{0 < \g \leq T} T^{i\ell\gamma}  (\Re{P}(\g))^k
     \frac{\mathop{d\omega}}{\omega} 
    +O(T)+O(T^{\frac12+\epsilon \ell}).
    \end{split}
    \ee
The  remaining  term here is handled in the same way as the corresponding term in \eqref{F Omega}, except that we use Conjecture~\ref{Conjecture} rather than Lemma~\ref{Land-Gon} to estimate the sums over $\g$. We carry this out now.
 
 Similarly to the analysis of the sum $S(k)$ in \eqref{S(k)}  that gave
 \eqref{contribution for k}, we find that
    \be\notag
    \begin{split}
    \sum_{0 < \g \leq T} T^{i\ell\gamma}  (\Re{P}(\g))^k
    =&\frac{1}{2^k}\sum_{j=0}^k \binom{k}{j} 
    \sum_{0 < \g \leq T} T^{i\ell\gamma}  \sum_n \frac{a_j(n)}{n^{1/2+i\g}} \sum_m \frac{a_{k-j}(m)}{m^{1/2-i\g}} \\
    =&\frac{1}{2^k}\sum_{j=0}^k \binom{k}{j}  \sum_n \frac{a_j(n)}{\sqrt n} \sum_m \frac{a_{k-j}(m)}{\sqrt m}
    \sum_{0 < \g \leq T} \Big(\frac{m T^\ell}{n}\Big)^{i\g},
    \end{split}
    \ee
where, as before,  $a_{r}(p_1\dots p_{r})$ denotes the number of permutations of the primes $p_1,\dots, p_{r}$, which might or might not be distinct.
By Conjecture~\ref{Conjecture},  for each $m$ and $n$
\[
\sum_{0 < \g \leq T} \Big(\frac{m T^\ell}{n}\Big)^{i\g} 
\ll T^{1+\ell \epsilon-\ell/2}  \Big(\frac{m}{n}\Big)^{-\frac12+\epsilon}
+T^{\frac12+\ell \epsilon}  \Big(\frac{m}{n}\Big)^{\epsilon}.
\]
Thus
\be \label{two sums}
\begin{split}
 & \sum_{0 < \g \leq T} T^{i\ell\gamma}  (\Re{P}(\g))^k \\
 \ll \,  & \frac{T^{1+\ell \epsilon-\frac{\ell}{2}}}{2^k}
 \sum_{j=0}^k \binom{k}{j}  \sum_n \frac{a_j(n)}{n^{\epsilon}}  
 \sum_m \frac{a_{k-j}(m)}{m^{1-\epsilon}}
 +
  \frac{T^{\frac12+\ell \epsilon} }{2^k}
  \sum_{j=0}^k \binom{k}{j} 
   \sum_n \frac{a_j(n)}{n^{1/2+\epsilon}} \sum_m \frac{a_{k-j}(m)}{m^{1/2-\epsilon}}\\
   &= \mathcal T_1(k)+ \mathcal T_2(k),
\end{split}
\ee
say.
Now 
\be \notag
\begin{split}
\mathcal T_1(k)  
&\ll  \frac{T^{1+\ell \epsilon-\frac{\ell}{2}}}{2^k}
 \sum_{j=0}^k \binom{k}{j}  \sum_n  {a_j(n)}  \sum_m {a_{k-j}(m)} \\
& \ll   \frac{T^{1+\ell \epsilon-\frac{\ell}{2}}}{2^k}\sum_{j=0}^k \binom{k}{j} 
    \Big( \sum_{p\leq X^2}1\Big)^j  \Big( \sum_{p\leq X^2} 1\Big)^{k-j}  
         \ll   T^{\frac12+\ell \epsilon }\pi(X^2)^k,
  \end{split}
 \ee
since $\ell \geq 1$.
Similarly,
\be \notag
\begin{split}
\mathcal T_2(k)  
&\ll   \frac{T^{\frac12+\ell \epsilon} }{2^k}
  \sum_{j=0}^k \binom{k}{j} 
   \sum_n  {a_j(n)}  \sum_m  {a_{k-j}(m)} 
   = T^{\frac12+\ell \epsilon} \pi(X^2)^k .
  \end{split}
 \ee

  Combining our estimates for $\mathcal T_1(k)$ and $\mathcal T_2(k)$ in \eqref{two sums}, we see that 
\[
 \sum_{0 < \g \leq T} T^{i\ell\gamma}  (\Re{P}(\g))^k
 \ll   T^{\frac12+\ell \epsilon}\  \pi(X^2)^k
\ll T^{\frac12+\ell \epsilon} X^{2k}  .
\]
Using this in \eqref{F Omega 3}, we find that
\be \notag
\begin{split}
     \sum_{0 < \g \leq T}  T^{i\ell\gamma} F_{\Omega}(\Re{P}(\g)-A) 
    &\ll  T^{\frac12+\ell \epsilon} \   X^{2K} \
    \int_0^\Omega  G\Big(\frac \omega \Omega\Big) 
    \sum_{1\leq k < K}  \frac{(2\pi\omega)^k}{k!} 
     \frac{\mathop{d\omega}}{\omega} 
     + O(T) +O(T^{\frac12+\epsilon \ell})\\
       &\ll    \Omega e^{2\pi \Omega} \ T^{\frac12+\ell \epsilon} \  X^{2K} 
     + O(T) .
\end{split}
\ee  
By \eqref{X, Omega} and \eqref{K} this is 
$$
\ll e^{7(\log\log T)^2} \ T^{\frac12+\epsilon \ell} \ T^{\frac{5}{(\log\log T)^{14}}}+T \ll T ,
$$
provided $\ell\leq L$ and $\epsilon$ is small enough relative to $L$.
Inserting this (and the same estimate when $A$ is replaced by $B$) into \eqref{set up 2}, we find that
\be\notag
\begin{split}
\sum_{0<\g \leq T}  
T^{i\ell  \g } \mathbbm{1}_{[A, B]}  &( \Re{P}(\g)) 
\ll T+ \frac{N(T)}{\Omega} \ll \frac{N(T)}{(\log\log T)^2} .
\end{split}
\ee
By \eqref{Weyl sum 1} and  \eqref{Weyl sum 2}  we then obtain
\be\notag
 \sum_{0<\g\leq T, \g\in \G_{[a, b]}} e\Big(\ell \g \frac{\log T}{2\pi}\Big)
 \ll N(T) \frac{(\log\log\log T)^2}{\sqrt{\log\log T}} 
\ee
for $\ell\leq L$.
Using  this bound in \eqref{ErdosTuran}, we  see that
\begin{equation}\notag
\Bigg| \sum_{\substack{0<\g \leq T, \g\in \G_{[a, b]}\\ \{ \g (\log T)/2\pi \} \in [\alpha,\beta]}} 1 \ \ - \ (\beta-\alpha)N_{[a, b]}(T) \Bigg| 
\ \ll \ \frac{N_{[a, b]}(T) }{L+1} \ + \   N(T) ( \log L) \frac{(\log\log\log T)^2}{\sqrt{\log\log T}}.
\end{equation} 
By our hypotheses on $a$ and $b$, $N_{[a,b]}(T)\gg N(T)$. Hence, since $L$ may be arbitrarily large, we find that this equals $o(N_{[a, b]}(T))$.
 This proves \eqref{discrp 1}. The analogous inequality when  $\G_{[a, b]}$ is replaced by $\G_{[a, b]}^*$ follows from this on noting that  
 the number of $\g$ in $\G_{[a, b]}$ with $0<\g\leq T$ that are not elements of $\G_{[a, b]}^*$
is at most $o(N_{[a,b]}(T))$.



\begin{thebibliography}{BM}

\normalsize
\baselineskip=16pt

 \bibitem{Baker} Baker, A. : \emph{Transcendental Number Theory}. Second edition. Cambridge Mathematical Library. Cambridge University Press, Cambridge, 1990. 

 \bibitem{log zeta}  \c{C}i\c{c}ek, F. : \emph{On the logarithm of the Riemann zeta-function near the nontrivial zeros}, Trans. Amer. Math. Soc., vol. 374, 2021, no. 8, 5995--6037.

 
\bibitem{Elliott} Elliott, P. D. T. A. : \emph{The Riemann zeta function and coin tossing}, J. Reine Angew. Math., vol. 254, 1972, 100--109.


%
%

\bibitem{GonekLandaulemma1} Gonek, S. M. : \emph{A formula of Landau and mean values of $\zeta(s)$}, Topics in analytic number theory (Austin, Tex., 1982), Univ. Texas Press, Austin, TX.(1985) 92--97.

\bibitem{GonekLandaulemma2} Gonek, S. M. : \emph{An explicit formula of Landau and its applications to the theory of the zeta-function}, Contemporary Mathematics 143, Amer. Math. Soc., Providence, RI, (1993) 395--413.





\bibitem{Mont} Montgomery, H. L. : \textit{Ten Lectures on the Interface between Analytic Number Theory and Harmonic Analysis}, CBMS Regional Conference Series in Mathematics, 84. Published for the Conference Board of the Mathematical Sciences, Washington, DC; by the American Mathematical Society, Providence, RI, 1994.




\bibitem{Rademacher} Rademacher, H. A. : \emph{Fourier Analysis in Number Theory}, 
Symposium on harmonic analysis and related integral transforms: Final technical report, Cornell Univ., Ithaca, N.Y. (1956).
 
 
\bibitem{T} Titchmarsh, E. C. : \emph{The Theory of the Riemann Zeta-function}, Oxford University Press, Second Edition, Oxford, 1986.

\bibitem{Tsang} Tsang, K. : \emph{The Distribution of the Values of the Riemann Zeta-function}. Thesis (Ph.D.)--Princeton University. 1984.  
 
\bibitem{Weyl} Weyl, H. : \emph{\"{U}ber die Gleichverteilung von Zahlen mod. Eins}, Math. Ann. 77 (1916), no. 3, 313--352.  
 
 
 
\end{thebibliography}
\end{document}